\magnification1200
\def\emptyset{\varnothing}
\mathchardef\[="505B
\mathchardef\]="405D

\font\tenBbb=msbm10  
 \textfont8=\tenBbb  
\def\Bbb{\fam=8}  
\def\R{{\Bbb R}}

\font\smallcaps=cmcsc10

\def\det{\mathop{\rm det}}   
 
\def\qed{\vbox{\hrule\hbox to7.8pt{\vrule height7pt   \hss\vrule
height7pt}\hrule} }

\font\bfuno=cmbx12
\font\smallcaps=cmcsc10

\def\Isom{\mathop{\rm Isom}}

\def\A{{\cal A}}

\mathchardef\emptyset="083F
\def\qed{%
  \vbox{\hrule\hbox to7.8pt{\vrule height7pt
  \hss\vrule height7pt}\hrule} }

\font\titolo=cmbx12

\centerline{{\titolo  DIFFEOMORPHISMS WITH}}
\centerline{{\titolo  BANACH SPACE DOMAINS}}

\bigskip                                                      
{\baselineskip14truept
\centerline{{\smallcaps Gaetano Zampieri}}
\centerline{{\it Universit\`a di Padova}}
\centerline{{\it Dipartimento  di Matematica
        Pura e Applicata}}
\centerline{{\it
              via   Belzoni 7, 
              35131 Padova, Italy}} }

Auxiliary scalar   coercive functions.
Sufficient
conditions for $f$ to be one-to-one. Sufficient conditions for 
bijectivity.

\bigskip\bigskip\bigskip
 
 \centerline{{\bfuno 1. Introduction}} 

 \bigskip

 The basic element of the following arguments is a $C^1$
 mapping  $f:X\to Y$, with $X,Y$ Banach spaces, and with  
 derivative   
 everywhere invertible: 
 $$f'(x)\in \Isom(X;Y)\qquad
 \forall x\in X\,.\eqno(1.1)$$
 So $f$ is a {\it local diffeomorphism} at every point by the
 {\it Inverse Function Theorem}.
 
 The aim of this paper is to find a sufficient condition for $f$
 to be injective, and so a {\it global diffeomorphism} $X\to f(X)$
 (Theorem~2.1), and a sufficient condition 
 for $f$ to be bijective and so a
 {\it global   diffeomorphism onto} $Y$ (Theorem~3.1). 
 This last condition is also necessary in the particular case
 $X=Y=\R^n$.
 
 In  these theorems the  key role is played by nonnegative {\it
 auxiliary
 scalar coercive functions}, that is, continuous mappings 
 $\;k:X\to\R_+$ with 
 $k(x)\to +\infty$ as $\|x\|\to +\infty\,$.
 As far as I
 know the use of such auxiliary functions in these questions is
 new. We find some first corollaries.  The author hopes that
 suitable auxiliary functions, adapted
 to particular problems, may lead to new  consequences.
 
 In order to briefly discuss the results, and some of their
 relations with the literature, let us consider the case
 where $X$ is a Hilbert space with scalar product `$\cdot$',
 and $k\in C^1$.
 This simplifies  the formulas a little.  However, our results
 will be formulated and proved for
 Banach spaces, where we will ask for $k$ to be locally
 Lipschitz continuous and to have 
 the right directional derivatives only; these last assumptions are
 not
 made to quibble about this matter, but  are
 related
 to the nondifferentiability of the map $x\mapsto \|x\|^2$ in
 general Banach spaces.
 
 In this paper we use the hypothesis that
 the operator norm $\|f'(x)^{-1}\|$ is bounded on bounded sets,
 i.e.
 $$\sup_{\|x\|\le r}\; \|f'(x)^{-1}\|<+\infty\,,\qquad
 \forall r:0<r<+\infty\,,\eqno(1.2)$$
 
 Theorem~2.1 (in the particular case we mentioned above)
 says that the local diffeomorphism $f$ is injective
 if (1.2) holds and there exist a point $x_0\in X$, and a coercive
 function $\;k\in
 C^1(X;\R_+)$, such that
 $$\sup_{x\in X}\ 
 k'(x)\,F(x)<+\infty\,,\qquad\hbox{{\rm with}}\quad
 F(x):=-f'(x)^{-1}(f(x)-f(x_0)).\eqno(1.3)$$
 
 By a suitable choice of $k$, 
 the following result  is achieved (Corollary~2.2): the local
 diffeomorphism $f$ is injective if (1.2) holds and
 there exist points $\,x_0, x_1\in X$, and non\-negative real
 numbers
 $a,b,c$, such that
 $$ (x-x_1)\cdot F(x)\le
 a+b\|x-x_1\|^2+c\|f(x)-f(x_0)\|^2\,,\qquad \forall x\in X
 \,.\eqno(1.4)$$ 
 
 To prove Theorem~2.1, the condition in (1.3) is used in connection
 with the following auxiliary o.d.e. (where $F$ is as in (1.3))
 $$\dot x=F(x)\,.\eqno(1.5)$$

 The equation~(1.5) plays an important role also
 in  [Zampieri, 1990] to prove the following sufficient
 condition of invertibility: the restriction of the local
 diffeomorphism
 $f:\Omega=\Omega^{\circ}\subseteq \R^n \to \R^n$ to the ball
 $\{x:|x-x_0|\le r\}\subset \Omega$, is one-to-one if 
 $$(x-x_0)\cdot F(x)\le 0\,, \qquad \forall x:
 |x-x_0|=r\,.\eqno(1.6)$$
 More generally, that paper suggests to estimate regions
 contained in the {\it `basin of attraction'} of $x_0$ for (1.5)
 to obtain domains of invertibility of $f$ around $x_0$. The
 results
 of [Zampieri, 1990] are generalized to Banach spaces by 
 Gianluca Gorni in [Gorni, 1990].
 
 In Section~3 we turn our attention to a  property stronger than
 one-to-oneness, namely bijectivity.  Several authors have dealt
 with bijectivity of local diffeomorphisms. We
 refer the reader
 to [Berger, 1977], [Ortega \& Rheinboldt, 1970], [Plastock, 1974],
[Prodi \& Ambrosetti,
 1973], and to [Ra\-du\-lescu
 \& Ra\-du\-lescu, 1980] for clear discussions on these topics, and
 some applications to Differential Equations.
 
 The  auxiliary functions $k$ are
 the `common denominator' of the present paper since they are also
 used 
 in Section~3 where we shall find a
 sufficient condition (Theorem~3.1) for bijectivity of local
 diffeomorphisms between Banach spaces.
 Theorem~3.1 says that the local diffeomorphism $f:X\to Y$ is 
 a global diffeomorphism onto $Y$ if (1.2) holds and there exists
 a coercive function $\,k\in C^1(X;\R_+)$ such that
 $$\sup_{x\in X}\ \|k'(x)\circ
 f'(x)^{-1}\|<+\infty\,.\eqno(1.7)$$ 
 In $\R^n$ the `if' becomes `if and only if' (and 
 condition (1.2) is always satisfied).
 As before, our result actually says much more, and we 
 now present the case $k\in C^1$ just for the sake of simplicity.

 From  Theorem~3.1 we first obtain a new  proof of the
 celebrated Theorem of Hadamard 
 generalized by
 Levy, Meyer, and Plastock to Banach spaces (Corollary~3.2
 below, and [Berger, 1977] Section~5.1).   This Theorem is based
 on a condition
 on  $\|f'(x)^{-1}\|$ which is assumed bounded or of slow growth
 in $\|x\|$ (roughly at most linear).

 Moreover, from Theorem~3.1 we  deduce another  known result:  
 $f$  is a  surjective
 diffeomorphism  if (1.1) and (1.2) hold, and 
 $f$ is coercive, i.e.
 $$\|f(x)\|\to +\infty\qquad\hbox{{\rm as}}\qquad \|x\|\to
 +\infty\,\eqno(1.8)$$
 (see Corollary~3.3 below and Theorem~2.1 in [Plastock,1974] which
 has essentially the same statement).
 
 As is well known, the local homeomorphism $f:X\to Y$ is bijective
 if and
 only if it is {\it proper}, i.e. $f^{-1}(C)$ is compact for any
 compact $C$. This last sentence is   an important Theorem
 (see [Berger, 1977] Theorem~5.14)
 which, for $X=Y=\R^n$, was glimpsed by Hadamard 
 (see [Hadamard, 1906] Section~17, and [Hadamard, 1968]).
 Later it was clarified and generalized    in
[Caccioppoli, 1932],
and in [Banach \& Mazur, 1934].

 In the particular case of local diffeomorphisms $\R^n\to \R^n$,
 the  framework of our Section~3, represents a unified point of
 view on
 the two celebrated  Theorems we mentioned above, the former
 being based on a condition on the growth of 
  $\|f'(x)^{-1}\|$ (assumed roughly at most linear), and the
latter based on properness of $f$, which in this case is equivalent
to coerciveness. Our Section~4 is devoted to the important case of
$\R^n$.

 The
mapping  $x\mapsto \arctan x$, $\;x\in \R$, can be proved to be
one-to-one  by Corollary~2.2 and it is not surjective.
A less trivial example  $f:\R^2\to
\R^2$ is shown in
Section 5. The criterion   
which includes (1.6) cannot be applied to prove the
injectivity of that $f$.

 Hopefully, suitable auxiliary functions can lead to new results
 of global bijectivity, as well as of plain injectivity, in
 particular in the realm of Differential Equations.
 \bigskip\bigskip\bigskip\goodbreak
 
 \centerline{{\bfuno 2. Proving injectivity of local
 diffeomorphisms}} 

 \bigskip
  
 We are going to use 
 nonnegative 
 auxiliary
 scalar coercive functions, that is continuous mappings 
 $\;k:X\to\R_+$ with 
 $k(x)\to +\infty$ as $\|x\|\to +\infty\,$.
 We also need the
 directional right derivatives, i.e. for any $x,v\in X$ we 
 assume the existence of
 $\ D^+_v k(x):=\lim_{s\to 0+}\;
 {k(x+sv)-k(x)\over {s}}\,$.
 Moreover, we need $k$ to be locally Lipschitz continuous.
 
 \bigskip
 \item{} \noindent {\bf Theorem~2.1}. $\ $\sl The mapping
 $f\in C^1(X;Y)$ is a global 
 diffeomorphism $X\to f(X)$ if (i) 
 $f'(x)\in \Isom(X;Y)\ \;
 \forall x\in X\,$, (ii)
 $$\sup_{\|x\|\le r}\; \|f'(x)^{-1}\|<+\infty\,,\qquad
 \forall r:0<r<+\infty\,,\eqno(2.1)$$
 and (iii) there exist a point $\,x_0\in X$ and
 a locally Lipschitzian coercive function $\, k:X\to \R_+$,
 which admits
 all $D^+_v k(x)$,  such
 that
 $$\sup\bigl\{ 
 D^+_v\,k(x):  v=-f'(x)^{-1}(f(x)-f(x_0)),\quad  x\in X 
 \bigr\}<+\infty\,.
 \eqno(2.2)$$
 Moreover, under these hypotheses, the range $f(X)$ is star shaped
 around $f(x_0)$, i.e.
 $$y\in f(X)\,,\ s\in [0,1]\qquad\Longrightarrow\qquad
 f(x_0)+s(y-f(x_0))\in f(X)\,.\eqno(2.3)$$
 For  $\,k\in C^1$  formula (2.2)
 becomes:
 $$\sup_{x\in X}\ 
 \bigl(-k'(x)f'(x)^{-1}(f(x)-f(x_0))\bigr)<+\infty\,.\eqno(2.4)$$
 
 \bigskip\rm

 {\bf Proof.}$\ $  Consider the Cauchy problem
 $$\dot x(t)=F(x(t)),\qquad x(0)=\bar x\,,\qquad \hbox{{\rm
with}}
 \quad F(x):=-f'(x)^{-1}(f(x)-f(x_0))\,.\eqno(2.5)$$
 
 We have that:
 
 \item{(a)}{there is local
 existence, uniqueness, and
 continuous dependence on the initial conditions for (2.5),} 
 \item{(b)}{if
 $\A$ is the `basin of attraction' of $x_0$, i.e. the set of all
 the
 points $\bar x\in X$ such that the
 maximal
 solution $t\mapsto x(t,\bar x)$ (with $x(0,\bar x)=\bar x$) is
 defined for all $t\ge 0$ and
 $x(t,\bar x)\to x_0$ as $t\to +\infty$, then 
 $\A$ is open and
 the restriction $f|\A$ is
 one-to-one, and} 
 \item{(c)}{if $\bar x\in \partial \A$, the boundary of
 $\A$, then the maximal solution $t\mapsto
 x(t,\bar x)$ cannot be global in the future, i.e. it cannot 
 be defined for all positive values of $t$.}
 
 A detailed proof of these facts can be found in
 the Lemmas of  [Gorni, 1990];
 here let us just
 remind the following fundamental facts. 
 We can check at once that, for any $(t,\bar x)$,
 $$ f\bigl(x(t,\bar x)\bigr)-f(x_0)\,=\,e^{-t}\,\bigl(f(\bar
 x)-f(x_0)\bigr)\,.\eqno(2.6)$$

 This  permits to prove (a) and (b). To have (c) we use
 (b) to say that $f|\A:\A\to f(\A)$ is a homeomorphism. This
 implies the existence of $d>0$ such that
 $\|f(x)-f(x_0)\|\ge d$ for all $x\in \partial \A$. Finally,
 (2.6) gives $t\le \ln (\|f(\bar x)-f(x_0)\|/d)$ for $\bar
 x\in \partial \A$ and $t$ in the domain of definition of 
 the maximal solution $t\mapsto x(t,\bar x)$. 
 
 Now, let us show that, under our hypotheses, the
 solutions to (2.5) are all global in the future. This implies that
 the boundary $\partial \A$ is empty (see (c) above). Thus $\A=X$
 and $f$ is one-to-one (see (b) above).
 
 Assume that  
 $t\mapsto x(t)$ is a solution to (2.5) (so in particular
 $\; x(0)=\bar x$) which is
 defined on $[0,b)$,
 with
 $0<b<+\infty$. We are going to prove that it can be extended to
 $\R_+$.

 By the hypothesis (iii) there exists a (nonnegative) function
 $k$ which satifies (2.2). 
 Let us define $\alpha: [0,b)\to \R_+,
 t\mapsto k(x(t))$.  Since (by hypothesis) $k$ admits all
 $D^+_v k(x)$
 and it is locally Lipschitzian, and 
 $t\mapsto x(t)\in C^1$, then the map $\alpha$ admits the
 right derivatives and we have
 $D^+\alpha(t)=D^+_{\dot x(t)}\,k(x(t))$. This is checked at
 once  by
 showing in particular that
 $$\lim_{s\to 0+}\;{k(x(t+s))-k(x(t)+s\dot x(t))\over s}=0\,.
 \eqno(2.7)$$
 Formulas (2.2) and (2.5) give
 $D^+\alpha(t)\le
 c$ where $c>0$ is the absolute value of the `$\sup$' in (2.2).
 Therefore, by a standard argument that we show below, we have that
 $$0\le k(x(t))=:\alpha(t)\le \alpha(0)+bc\,.\eqno(2.8)$$
 So $x(t)\in k^{-1}([0,\alpha(0)+bc])$. This last set is contained
 in
 a ball,
 say $\|x\|\le \hat r$, since $k:X\to \R_+$ is coercive ($k(x)\to
 +\infty$ as $\|x\|\to +\infty$). Thus
 $$\|x(t)\|\le \hat r\qquad \forall t\in
 [0,b)\,.\eqno(2.9)$$
 One of our hypotheses is (2.1), so we can define
 $ a:=\sup_{\|x\|\le \hat r}\; \|f'(x)^{-1}\|<+\infty\,$. Moreover
 the map $t\mapsto \|f(x(t))-f(x_0)\|$  is decreasing (see (2.6)).
 From (2.9) and (2.5)  we then have 
 $\|\dot
 x(t)\|\le a\,\|f(\bar x)-f(x_0)\|\,$.
 Thus $t\mapsto x(t)$ is Lipschitzian, and   has a limit as
 $t\to b$. So it
 can be extended in the future, and 
the maximal solution is defined for all positive
$t$.

 In the particular case where $k\in C^1$, formula (2.2) gives
 (2.4)
 at once. 
 
 Finally, let $\bar y\in f(X)$, and $s\in (0,1]$. We define $\bar
 x=f^{-1}(\bar y)$ (we just proved that $f$ is 1-1), and consider
 the maximal solution $t\mapsto x(t,\bar x)$. Since $X=\A$, we
 have that this is globally defined in the future. So we can
 consider $t:=\ln(1/s)$ and (2.6) proves (2.3) for $s\ne 0$. The
 case $s=0$ is trivial.
 
 All we are left to prove is formula (2.8) from 
 $D^+\alpha(t)\le c$.
 We fix any $\epsilon
 >0$. It is enough to 
 prove that $\alpha(t)-\alpha(0)\le (c+\epsilon)t$.
 The set where this holds is an interval $I_{\epsilon}$. We argue
 by contradiction and assume that
 $\hat b:=\sup I_{\epsilon}<b$. Then, by the continuity of
 $\alpha$, we have that $\hat b\in I_{\epsilon}$. The existence
 of $D^+\alpha(\hat b)$ implies that
 $${\alpha(t)-\alpha(\hat b)\over{t-\hat b}}\le D^+\alpha(\hat
 b)+\epsilon\le  c+\epsilon$$
 for $\hat b<t<b$, with $t$
 near $\hat b$. So,  for such values of $t$, we have  (remind that
 $\hat b\in I_{\epsilon}$)\hfil\break
 $\alpha(t)-\alpha(0)=(\alpha(t)-\alpha(\hat b))+(\alpha(\hat
 b)-\alpha(0))\le (c+\epsilon)(t-\hat b)+ (
 c+\epsilon)\hat b=
 (c+\epsilon)t$. 

 \line{\hfil\qed}
 \bigskip
 \goodbreak
 
 \item{} \noindent {\bf Corollary~2.2}. $\ $\sl The mapping
 $f\in C^1(X;Y)$ is a global 
 diffeomorphism $\ X\to f(X)$ if  conditions (i) and (ii) in
 Theorem~2.1 are satisfied, and
 there exist points $\,x_0, x_1\in X$, and nonnegative real numbers
 $a,b,c$,
 such that
 $$\eqalign{&D^+_{F(x)}\,g(x)\le
 a+b\|x-x_1\|^2+c\|f(x)-f(x_0)\|^2\,,\qquad\quad\forall x\in
 X\,,\cr
 &\hbox{{\sl
 where}}\quad F(x):=-f'(x)^{-1}(f(x)-f(x_0)),\quad \hbox{{\sl
 and}}\quad
 g(x):=\|x-x_1\|^2\,.\cr} 
 \eqno(2.10)$$
 Moreover (2.3) holds. If $X$ is a Hilbert space then
 formula (2.10) becomes:
 $$ -2(x-x_1)\cdot f'(x)^{-1}(f(x)-f(x_0))\le
 a+b \|x-x_1\|^2+c\|f(x)-f(x_0)\|^2\,,\quad \forall x
 \eqno(2.11)$$

 \bigskip\rm
 \goodbreak

 {\bf Proof.}$\ $  We can assume that $a\ge b>2$ (if this is false,
 we may use
 $a+b+3$ and $b+3$ instead of $a$ and $b$ respectively). Consider
 the
 following auxiliary function
 $$k(x):=\ln h(x)\qquad\hbox{{\rm with}}\qquad
 h(x):={a\over b}
 +\|x-x_1\|^2+{c\over {b-2}}\|f(x)-f(x_0)\|^2\,.\eqno(2.12)$$
 This is trivially coercive. Moreover, 
 it is locally Lipschitzian and it admits all $D^+_vk(x)$. We have
 $$\eqalign{h(x)\ D^+_{F(x)}&\,k(x)=\;
 D^+_{F(x)}\,g(x)+\cr
 &+2 {c\over {b-2}} \|f(x)-f(x_0)\|\  \lim_{s\to
 0+}{\|f(x+sF(x))-f(x_0)\|-\|f(x)-f(x_0)\|\over{s}}\,,\cr}$$
 where $D^+_{F(x)}\,g(x)$ and the other limit exist as one checks 
 by using the convexity of the norm. So (2.10) gives
 $$\eqalign{h(x)\ D^+_{F(x)}&\,k(x)\le
 a+b\|x-x_1\|^2+c\|f(x)-f(x_0)\|^2+\cr
 &+2 {c\over {b-2}}\|f(x)-f(x_0)\|\; \|f'(x) F(x)\|=b h(x)\,.\cr}$$
 This shows that (2.2) holds. Finally Theorem~2.1 gives 
 Corollary~2.2.
 
 \line{\hfil\qed}
\bigskip\bigskip\bigskip\goodbreak


 \centerline{{\bfuno 3. Proving bijectivity of local
 diffeomorphisms}}

 \bigskip
 \item{} \noindent {\bf Theorem~3.1}.$\ $\sl The mapping
 $f\in C^1(X;Y)$ is a global 
 diffeomorphism onto $Y$ if (i) 
 $f'(x)\in \Isom(X;Y)\ $ $\;
 \forall x\in X\,$, (ii)
 $$\sup_{\|x\|\le r}\; \|f'(x)^{-1}\|<+\infty\,,\qquad
 \forall r:0<r<+\infty\,,\eqno(3.1)$$
 and (iii) there exists 
 a  locally Lipschitzian coercive function $\, k:X\to \R_+$
 which admits
 all $D^+_v k(x)$, and it is such
 that
 $$\sup\bigl\{ 
 D^+_v\,k(x):  v=f'(x)^{-1}u,\  x\in X,\; u\in Y,\;
 \|u\|=1\bigr\}<+\infty\,. \eqno(3.2)$$
 For  $\,k\in C^1$ this last formula 
 is equivalent to:\quad
 $\sup_{x\in X}\ 
 \|k'(x)\circ f'(x)^{-1}\|<+\infty\,$. 
 \bigskip\rm

 \goodbreak
 {\bf Proof.$\ $}  
 Consider the Cauchy problem (2.5), that is
  $$\dot x(t)=F(x(t)),\qquad x(0)=\bar x\,,\qquad \hbox{{\rm
with}}
 \quad F(x):=-f'(x)^{-1}(f(x)-f(x_0))\,\eqno(3.3)$$
 where $\bar x, x_0$ are
any
 distinct points. Its maximal solution $t\mapsto x(t,\bar x)$
satisfies  (2.6), i.e. 
 $$ f\bigl(x(t,\bar x)\bigr)-f(x_0)\,=\,e^{-t}\,\bigl(f(\bar
 x)-f(x_0)\bigr)\,.\eqno(3.4)$$
 Since 
 $$t\mapsto \|f\bigl(x(t,\bar x)\bigr)-f(x_0)\|$$
 is bounded whenever $t$ ranges on a bounded interval, we may just
 repeat some arguments of the proof of Theorem~2.1 to have the
global existence in the future of the solution to our Cauchy
problem. In these arguments we consider the derivative
$D_v^+\,k(x)$  with 
$$v=-f'(x)^{-1}u\,,\qquad u={f(\bar x)-f(x_0)
\over{\|f(\bar x)-f(x_0)\|}}$$
(see (3.3) and remind (3.4)). The global existence in the future
implies the injectivity of $f$ as we saw in the proof of
Theorem~2.1.
 
 Our actual hypothesis (3.2), unlike the one of Theorem~2.1, permits
to say that the solution to the Cauchy problem (3.3) is global in
the past too. Indeed, we just need to consider the opposite vector
field,  namely the differential equation
$$ \dot z(t)=-F(z(t))\,,\qquad \quad \hbox{{\rm with}}\quad
F(x)=-f'(x)^{-1}(f(x)-f(x_0))\,,\eqno(3.5)$$
whose solution starting at $\bar x$ is $z(t)=x(-t,\bar x)$. Now
 the map 
$t\mapsto \|f(z(t))-f(x_0)\|$ increases, but what we need is
boundedness on bounded intervals only.

So the solution to (3.3) is defined on the whole $\R$ and $f$ maps
it to the half-line
$$\{y\in Y: y=f(x_0)+\xi \bigl(f(\bar x)-f(x_0)\bigr),\  \;
0<\xi<+\infty\}$$  (see (3.4)). This implies the 
 surjectivity of $f$ since $f(\bar x)$ ranges in a full
neighbourhood of $f(x_0)$ ($f(x_0)$ excluded).

 \line{\hfil\qed}
 \bigskip

 \goodbreak
 The following theorem is  known. We find it again as a
 consequence of Theorem~3.1.

 \bigskip
 \item{} \noindent {\bf Corollary~3.2}.$\ $\sl 
 Let
 $f\in C^1(X;Y)$, $f'(x)\in \Isom(X;Y)\ \;
 \forall x\in X\,$. Then
 $f$ is a global diffeomorphism onto $Y$ if 
 there
 exists a continuous  map
 $\omega:\R_+\to\R_+\setminus\{0\}$
 such that
 $$\int_0^{+\infty}\, {1\over{\omega(s)}}\,
 ds=+\infty\,,\qquad
 \|f'(x)^{-1}\|\le \omega(\|x\|)\,.\eqno(3.6)$$
 \item{}  In particular this holds if, for some
 $a,b\in \R_+\,$, we have
 $$\  \|f'(x)^{-1}\|\le a +b \|x\|\,.\eqno(3.7)$$ 
 
 \bigskip\rm
 
 \goodbreak
 {\bf Proof.$\ $}
 Define
 $$k:X\to\R_+,\  x\mapsto\int_0^{\|x\|}\,
 {1\over{\omega(s)}}\, ds\,.\eqno(3.8)$$
 By the first condition in (3.6)
 we have that this function  is coercive. Moreover it is 
 locally Lipschitz
 continuous and it admits all
 the derivatives $D^+_v k(x)$. For $v=f'(x)^{-1}u$, $x\in X$,
 $u\in Y$, $\|u\|=1$,  we have (see (3.6))
 $$D^+_v
 k(x)={1\over{\omega(\|x\|)}}\; \lim_{s\to 0+}
 {\|x+sv\|-\|x\|\over{s}}\;\le
 \;{1\over{\|f'(x)^{-1}\|}}\;\|v\|\;\le
 1\,.$$
 Where the limit exists as one verifies by means of the
 convexity of the norm. So (3.2) is satisfied, (3.1) holds by
 (3.6), and Theorem~3.1
 gives  Corollary~3.2. 
 
 \line{\hfil\qed}
  
 \goodbreak
  
 \bigskip \noindent

 As we said in Section~1, Corollary~3.2
 was discovered by Hadamard in
 $\R^n$.
 For  Banach spaces it was proved
 in [Levy, 1920] under condition (3.7) with $b=0$; 
 [Meyer, 1968] demonstrated that (3.7) is sufficient,
 and finally  [Plastock, 1974] gave a proof for the 
 general statement.
 
 Now, we are going to see another known consequence of 
 Theorem~3.1.

 \bigskip
 \item{} \noindent {\bf Corollary~3.3}.$\ $\sl The mapping
 $f\in C^1(X;Y)$ is a global surjective
 diffeomorphism if (i) it is coercive, i.e. 
 $\,\|f(x)\|\to +\infty\,$ as $\,\|x\|\to
 +\infty\,$, (ii)
 $\,f'(x)\in \Isom(X;Y)\ \;  \forall x\in X\,$,
 and (iii) the condition in (3.1) is satisfied.
 
 \goodbreak
 \bigskip\rm
 {\bf Proof.$\ $}  We define the function
 $$k(x):=\ln(1+\|f(x)\|^2)\,.\eqno(3.9)$$
 This is  coercive as well as $f$. Furthermore it admits all 
 $D^+_v k(x)$ and it is locally Lipschitz continuous. For
$v=f'(x)^{-1}u$,
 $x\in X$,
 $u\in Y$, $\|u\|=1$,  we have  
 $$D^+_v
 k(x)={2\|f(x)\|\over{1+\|f(x)\|^2}}\; \lim_{s\to 0+}
 {\|f(x+sv)\|-\|f(x)\|\over{s}}\le
 {2\|f(x)\|\over{1+\|f(x)\|^2}}\;\|f'(x)v\|
 \le 1\,.\eqno(3.10)$$
 So the condition (3.2) holds and Theorem~3.1 gives Corollary~3.3.
 
 \line{\hfil\qed}
\bigskip\bigskip\bigskip\goodbreak

 \centerline{{\bfuno 4. Finite dimension}}

 \bigskip
 
  Of course the case of the Euclidean space is particularly
  important and relevant  applications of global
  inverse function theorems in finite dimension  arise in Numerical
Analysis, Network
  Theory, Economics and other fields (see [Sandberg, 1980]).
  Let me also mention the Jacobian conjecture
  for global
  asymptotic stability. In the plane this conjecture leads to
  an injectivity problem which is still open (see [Zampieri \&
  Gorni, 1991]).
  
 In $\R^n$ all we have said  is simpler. 
 The
 auxiliary functions $k$ in
 Theorems~2.1 and 3.1 
 may be taken $C^1$. 
 Let us see how the previous results can be stated in $\R^n$.
 \bigskip
\goodbreak
\item{} {\bf Theorem~4.1}.\ \ \sl Let 
$$f:\R^n \to \R^n\,,\qquad
f\in C^1\,,\qquad \det f'(x)\ne 0\,,\qquad \forall x\in
\R^n\,.\eqno(4.1)$$
 Then $f$ is
one-to-one if there exists a coercive function $\ k\in
C^1(\R^n;\R_+)$, and a point $x_0\in \R^n$ such that
$$ \sup\ \{k'(x)\; F(x): x\in \R^n\}<+\infty ,
\quad \hbox{{\sl
with}}\quad 
F(x):=-f'(x)^{-1}\,\left(f(x)-f(x_0)\right).\eqno(4.2)$$
\bigskip\goodbreak  \rm \noindent
If one prefers $\ k'(x)\; F(x)=\nabla k(x)\cdot F(x)$, i.e.
the scalar
product of the gradient of $k$ and the vector field $F$. This
theorem still gives Corollary~2.2 in the finite dimension, i.e.
the following

\bigskip
\goodbreak
\item{} {\bf Corollary~4.2}.\ \ \sl Let $f$
be as in formula (4.1).  Then  $f$ is
one-to-one if there exist points $x_0, x_1 \in \R^n$, and
nonnegative real numbers $\; a,b,c,\;$ such that
$$ (x-x_1)\cdot F(x)\le
a+b |x-x_1|^2+c |f(x)-f(x_0)|^2\,,\qquad \forall x\in \R^n
\,,\eqno(4.3)$$ 
where $F(x)$ is as in (4.2).\rm
\bigskip\goodbreak  \rm 
Now let us turn our attention to  bijectivity. 
 The `if' in Theorem~3.1 and Corollary~3.3 can be substituted by
`if and
 only
 if'.   To verify the necessity we just remark that if $f:\R^n
 \to \R^n$ is a
 global surjective diffeomorphism then
 $f'(x)\in \Isom(\R^n;\R^n)\ \;
 \forall x\in \R^n\,$ and $f$ is coercive. 
 Moreover we can just consider the mapping in
 (3.9), which is coercive as well as $f$, and (3.10) completes
 the argument.
 
\bigskip
\goodbreak
\item{} {\bf Theorem~4.3}.\ \ \sl 
 The mapping $f\in C^1(\R^n;\R^n)$ is
a (global) diffeomorphism onto $\R^n$ if and only if (i)
$\det f'(x)\ne 0$ at every $x\in \R^n$, and (ii) there exists a
coercive function $\ k\in C^1(\R^n;\R_+)$ such that
$$ \sup\ \{\|k'(x)\circ f'(x)^{-1}\|: x\in \R^n\}<+\infty \,.
\eqno(4.4)$$
In other words, (4.4) may be written as
$$\sup\ \{|\nabla k(x)\cdot f'(x)^{-1}\,u|: x\in
\R^n, u\in \R^n, |u|=1\}<+\infty \,. \eqno(4.5)$$ 
\bigskip\goodbreak  \rm 
\noindent
From the last theorem we can easily deduce
the following two celebrated results.
The  proofs above can be easily adapted to $k\in C^1$ (of
 course
 the map in (3.8) is not differentiable at $x=0$ but we can
 just define $k$ in a different way for $|x|\le 1$).
 
\bigskip
\goodbreak
\item{} {\bf Corollary~4.4}.\ \ \sl   Let $f$
be as in formula (4.1).  Then  $f$ is
a  diffeomorphism onto $\R^n$ if 
there
 exists a continuous  function
 $\omega:\R_+\to\R_+\setminus\{0\}$
 such that
 $$\int_0^{+\infty}\, {1\over{\omega(s)}}\,
 ds=+\infty\,,\qquad
 \|f'(x)^{-1}\|\le \omega(|x|)\,.\eqno(4.6)$$
 \item{}  In particular this holds if, for some
 $a,b\in \R_+\,$, we have
 $$\  \|f'(x)^{-1}\|\le a +b |x|\,.\eqno(4.7)$$ 
 \bigskip\goodbreak  \rm 
 Finally we have the following known Theorem.
 \bigskip
\goodbreak
\item{} {\bf Corollary~4.5}.\ \ \sl  
The mapping $f\in C^1(\R^n;\R^n)$ is
a (global) diffeomorphism onto $\R^n$ if and only if (i)
$\det f'(x)\ne 0$ at every $x\in \R^n$, and (ii)
 it is coercive,
i.e.
$$ |f(x)|\to +\infty,\qquad \hbox{{\sl as}}\qquad |x|\to
+\infty.\eqno(4.8)$$
\bigskip\goodbreak  \rm 

 \bigskip\bigskip\bigskip\goodbreak


\centerline{{\bfuno 5. A nonsurjective example}}

 \bigskip

Let us give an example $\; f:\R^2\to \R^2$ which
satisfies the condition (4.3) in Corollary~4.2  but
which is not onto $\R^2$. Moreover, in this example
the left hand side of (4.3) assumes values of both
signs on every circumference $\;|x|=r>0$, so 
the criterion   
which includes (1.6) cannot be applied to prove the
injectivity.

In the sequel $x=(\xi,\eta)^T\in \R^2$ and our
function is
$$ f:\R^2\to \R^2\,,\quad\qquad {\xi\choose \eta}\;
\mapsto\;
{e^{\xi}\over{\sqrt{1+\eta^2}}}\,{1\choose
\eta}\,.\eqno(5.1)$$
This mapping is not surjective since 
its first component is positive.

We are going to prove that condition (4.3) for
injectivity is satisfied if we choose  $x_0=x_1$ at
the origin,  $\;a=b=1$, and $c=0$; namely
$$x\cdot F(x)\,\le\,
1+|x|^2\,,\qquad\quad\hbox{{\rm where}}\quad 
F(x)=-f'(x)^{-1}\,\bigl(f(x)-f(0)\bigr)\,.
\eqno(5.2)$$

We have
$$f'(x)=
{e^{\xi}\over{{(1+\eta^2})^{3\over 2}}}\;
\left(
  \matrix{1+\eta^2 & -\eta\cr
  \noalign{\medskip}
 \eta(1+\eta^2)& \ 1\cr}
  \right)\,,
  \eqno(5.3)$$
whose  determinant nowhere vanishes.
Moreover 
$$f'(x)^{-1}=
{e^{-\xi}\over{\sqrt{1+\eta^2}}}\;
\left(
  \matrix{1 & \eta\cr
  \noalign{\medskip}
 -\eta(1+\eta^2)& \ 1+\eta^2\cr}
  \right)\,,
  \eqno(5.4)$$
$$F(x)={e^{-\xi}\over{\sqrt{1+\eta^2}}}\;
{1\choose -\eta(1+\eta^2)}\,-\,{1\choose
0}\,.\eqno(5.5)$$
Therefore the left hand side of the inequality
 that
we are checking is
$$x\cdot F(x)=\xi\,
\left[{e^{-\xi}\over{\sqrt{1+\eta^2}}}-1\right]\,-\,
\eta^2 e^{-\xi}\sqrt{1+\eta^2}\,.\eqno(5.6)$$
If (i) $\xi\ge 0$ or (ii) $\xi<0$ and the term
between brackets in (5.6) is nonnegative, then 
$x\cdot F(x)\le 0$ and formula (5.2) holds.
Otherwise we have (iii) 
$\xi<0$ and the term between brackets is
negative too (and greater than $-1$). In this last
case we have $$x\cdot F(x)\,\le\, \xi\,
\left[{e^{-\xi}\over{\sqrt{1+\eta^2}}}-1\right]\,
\le\,
|\xi|\,\le\, 1+\xi^2\,\le\,
1+\xi^2+\eta^2\,.\eqno(5.7)$$

\bigskip\bigskip\bigskip\goodbreak


 \centerline{{\bfuno  Acknowledgements}} 
 \bigskip
The author thanks Giuseppe De Marco for a
critical discussion on the manuscript.
Furthermore, he thanks the {\it ``Ministero
dell'Universit\`a e della Ricerca Scientifica e
Tecnologica''}  which
supported this research.

\bigskip\bigskip\bigskip

 \centerline{{\bfuno  References}}  
 \bigskip

 \frenchspacing
 \item{1.} {\smallcaps Banach, S., \& Mazur, S.} 
 {\sl
 \"Uber mehrdeutige stetige Abbildungen}. {\it
 Studia Math. 5}, 174--178 (1934).
 
 \medskip
 
 \item{2.} {\smallcaps Berger, M.S.}   {\it 
 Nonlinearity and Functional Analysis}. Academic
 Press (1977).
 
 \medskip

 \item{3.} {\smallcaps Caccioppoli, R.}  {\sl
 Sugli
 elementi uniti delle trasformazioni funzionali:
 un teorema di esistenza e di unicit\`a ed alcune
 sue applicazioni}. {\it Rend. Sem. Mat.
 Padova~3}, 1--15 (1932).

 \medskip

 \item{4.} {\smallcaps Gorni, G.} {\sl A
 criterion of
 invertibility in the large for local
 diffeomorphisms between Banach spaces}. {\it Nonlinear Anal. 21} (1993), no. 1, 43Ð47. 
 \medskip

 \item{5.} {\smallcaps Hadamard, J.}  {\sl Sur
 les transformations ponctuelles}. {\it Bull.
 Soc. Math. France 34}, 71--84 (1906).
 
 \medskip

 \item{6.} {\smallcaps Hadamard, J.} {\sl Sur
 les correspondances ponctuelles}. {\it
 Oeuvres~I}, Editions du CNRS,
 383--384 (1968).
 
 \medskip

 \item{7.} {\smallcaps Levy, P.}  {\sl Sur les
 fonctions des lignes implicit\'es}. {\it Bull.
 Soc. Math. France~48}, 13--27 (1920).
 
 \medskip
 
 \item{8.} {\smallcaps Meyer, G.} {\sl On solving
 nonlinear equations with a one parameter
 imbedding}. {\it Siam J. Numer. Anal.~5},
 739--752 (1968).

 \medskip

\item{9.} {\smallcaps  Ortega, J.M. \& 
Rheinboldt, W.C.} {\it Iterative solutions of
nonlinear equations in several variables}.
Academic Press (1970).
 
 \medskip

 \item{10.} {\smallcaps Plastock, R.}  {\sl
 Homeomorphisms between Banach spaces}. {\it
 Trans. Am. Math Soc.~200}, 169--183 (1974).
 
 \medskip

 \item{11.} {\smallcaps Prodi, G., \& Ambrosetti,
 A.} {\it
 Analisi non lineare}. Quaderni della Scuola
 Normale Superiore, Pisa, Italy
 (1973).

 \medskip

 \item{12.} {\smallcaps Radulescu,
 M., \& Radulescu, S.}  {\sl Global
 inversion theorems and applications to
 differential equations}.
 {\it Nonlinear Analysis~4}, 951--965 (1980).
 
 \medskip

\item{13.} {\smallcaps  Sandberg, I. W.} {\sl
Global inverse function theorems}. {\it I.E.E.E.
Trans. Circuits Systems CAS~27}, No. 11,
 998--1004 (1980).

\medskip

\item{14.} {\smallcaps Zampieri, G., \& Gorni,
G.}  {\sl On the Jacobian conjecture for
global asymptotic stability}. {\it J. Dynam. Differential Equations 4} (1992), no. 1, 43Ð55.

\medskip
 
 \item{15.}  {\smallcaps Zampieri, G.}
 {\sl Finding domains of invertibility for
 smooth functions by means of attraction 
 basins}.  {\it J. Differential Equations 104} (1993), no. 1, 11Ð19. 

 \bye